\documentclass[a4paper, 12pt]{extarticle}
\usepackage{cite}
\usepackage{setspace}
\usepackage{mathtext}
\usepackage{graphicx}
\usepackage{amscd}
\usepackage[T2A]{fontenc}
\usepackage[english]{babel}
\usepackage{latexsym,amsfonts,amssymb,amsmath,longtable,amsthm}
\usepackage{cite}
\let\OLDthebibliography\thebibliography
\renewcommand\thebibliography[1]{
  \OLDthebibliography{#1}
  \setlength{\parskip}{0pt}
  \setlength{\itemsep}{0pt plus 0.3ex}
}
\usepackage[pic]{xy}
\usepackage{bbm,makeidx}
\usepackage[centerlast,small]{caption}
\usepackage{amsmath}
\usepackage[cp1251]{inputenc}
\usepackage[usenames, dvipsnames]{color}
\newtheorem{lemma}{{\scshape Lemma}}
\newtheorem{cor}{{\scshape Corollary}}

\newtheorem{theorem}{{\scshape Theorem}}

\begin{document}
\title{Chevalley groups of types $B_n$, $C_n$, $D_n$ over certain fields do not possess the $R_{\infty}$-property}
\author{Timur Nasybullov\footnote{The author is supported by the Research Foundation -- Flanders (FWO): postdoctoral grant  12G0317N and travel grant V425219N.}}
\date{}
\maketitle
\begin{abstract}
Let $F$ be an algebraically closed field of zero characteristic. If the transcendence degree of $F$ over $\mathbb{Q}$ is finite, then all Chevalley groups over $F$ are known to possess the $R_{\infty}$-property. If the transcendence degree of $F$ over $\mathbb{Q}$ is infinite, then Chevalley groups of type $A_n$ over $F$ do not possess the $R_{\infty}$-property. In the present paper we consider Chevalley groups of classical series $B_n$, $C_n$, $D_n$ over $F$ in the case when the transcendence degree of $F$ over $\mathbb{Q}$ is infinite, and prove that such groups do not possess the $R_{\infty}$-property. \\

\end{abstract}
\section{Introduction}
Let $G$ be a group and $\varphi$ be an automorphism of $G$. Elements $x$, $y$ from $G$ are said to be $\varphi$-conjugated if there exists an element $z\in G$ such that $x=zy\varphi(z)^{-1}$. The relation of $\varphi$-conjugation is an equivalence relation and it divides $G$ into  $\varphi$-conjugacy classes. The number $R(\varphi)$ of these classes is called the Reidemeister number of $\varphi$.  

Twisted conjugacy classes appear naturally in Nielsen-Reidemeister fixed point theory. Let $X$ be a finite polyhedron and $f:X\to X$ be a homeomorphism of $X$. Two fixed points $x,y$ of $f$ are said to belong to the same fixed point class of $f$ if there exists a path $c:[0,1]\to X$ with $c(0)=x$ and $c(1)=y$ such that $c\simeq f\circ c$, where $\simeq$ denotes the homotopy with fixed endpoints. The relation of being in the same fixed point class is an equivalence relation on the set of fixed points of $f$. The number $R(f)$ of fixed point classes of $f$ is called the Reidemeister number of the homeomorphism $f$. The Reidemeister number $R(f)$ is the homotopic invariant of $f$ and it plays a crucial role in the Nielsen-Reidemeister fixed point theory \cite{Kia}. Denote by $G=\pi_1(X)$ the fundamental group of $X$, and by $\varphi$ the automorphism of $G$ induced by $f$. In this notation the number $R(f)$ of fixed point classes of $f$ is equal to the number $R(\varphi)$ of $\varphi$-conjugacy classes in $G$ (see \cite[Chapter~III, Lemma~1.2]{Kia}). Thus, the topological problem
of finding $R(f)$ reduces to the purely algebraic problem of finding $R(\varphi)$. 

The Reidemeister number is either a positive integer or infinity and we do not distinguish different infinite cardinal numbers denoting all of them by the symbol~$\infty$. If $R(\varphi)=\infty$ for all automorphisms $\varphi$ of $G$, then  $G$ is said to  possess the $R_{\infty}$-property. The problem of determining groups which possess the $R_{\infty}$-property was formulated by A.~Fel'shtyn and R.~Hill \cite{FelHil}. The study of this problem
has been quite an active research topic in recent years. We refer to the paper \cite{FelNas} for an overview of the families of groups which have been studied in this context until 2016. More recent results can be found in \cite{DekGon2,MubSan, Nas5, Tro,Tro2,Yoo}. For the immediate consequences of the $R_{\infty}$-property for topological fixed point
theory see~\cite{GonWon}. Some aspects of the $R_{\infty}$-property can be found in \cite{FelTro}.

The author studied conditions which imply the $R_{\infty}$-property for different linear groups over rings \cite{Nas1, Nas3,  Nas5} and fields \cite{FelNas, Nas2, Nas4, Nas5}. 
In particular, it was proved that if $F$ is a field of zero characteristic which has either finite transcendence degree over $\mathbb{Q}$, or periodic group of automorphisms, then every Chevalley group (of normal type) over $F$ possesses the $R_{\infty}$-property \cite{Nas2}. In the present paper we consider Chevalley groups of classical series $B_n, C_n, D_n$ over algebraically closed fields of zero characteristic which have infinite transcendence degree over $\mathbb{Q}$, and prove the following main result.

~\\
\textbf{{\scshape Theorem \ref{main}.}} {\it Let $F$ be an algebraically closed field of zero characteristic which has infinite transcendence degree over $\mathbb{Q}$. Then Chevalley groups of types $B_n, C_n, D_n$ over $F$ do not possess the $R_{\infty}$-property.}\\

In \cite[Theorem~7]{Nas4} a similar result was proved for Chevalley groups of type $A_n$, therefore we do not consider the case of root systems $A_n$ in the present paper. In order to prove Theorem~\ref{main} we use classical matrix representations of adjoint Chevalley groups of types $B_n$, $C_n$, $D_n$. Chevalley groups of types $E_6$, $E_7$, $E_8$, $F_4$, $G_2$ are not isomorphic to certain classical linear groups and we cannot use the same technique.

 Theorem~\ref{main} together with the results \cite[Theorem~3.2]{FelNas} and \cite[Theorem~7]{Nas4} imply the following result which gives necessary and sufficient condition of the $R_{\infty}$-property for Chevalley groups of classical series $A_n, B_n, C_n, D_n$ over algebraically closed fields.

~\\
\textbf{{\scshape Theorem \ref{main2}.}} {\it Let $F$ be an algebraically closed field of zero characteristic. Then Chevalley groups of types $A_n,B_n, C_n, D_n$ over $F$ possess the $R_{\infty}$-property if and only if the transcendence degree of $F$ over $\mathbb{Q}$ is finite.}\\

The paper is organized as follows. In Section~\ref{sec1} we give necessary preliminaries from field theory and recall the construction from \cite{Nas4} of one specific automorphism of the field $\overline{\mathbb{Q}(s_1,s_2,\dots)}$, where $s_1,s_2,\dots$ is a countable set of variables. In Sections~\ref{sec3},~\ref{sec4} we prove that orthogonal and symplectic linear groups over $\overline{\mathbb{Q}(s_1,s_2,\dots)}$ do not possess the $R_{\infty}$-property. In Section~\ref{sec5} we prove that  orthogonal and symplectic linear groups over an arbitrary algebraically closed field $F$ of zero characteristic with infinite transcendence degree over $\mathbb{Q}$ do not possess the $R_{\infty}$-property. In Section~\ref{sec6} we consider classical Chevalley groups over an arbitrary algebraically closed field of zero characteristic with infinite transcendence degree over $\mathbb{Q}$ and prove Theorem~\ref{main} and Theorem~\ref{main2}.
\section{One automorphism of $\overline{\mathbb{Q}(s_1, s_2,\dots)}$}\label{sec1}

Let $L$ be a subfield of a field $F$, and $X$ be a subset of $F$. The minimal subfield of $F$ which contains $L$ and $X$ is denoted by $L(X)$. Elements $x_1,\dots,x_k$ of $F$ are called  algebraically independent over $L$ if there is no polynomial $f(t_1,\dots,t_k)\neq0$ with coefficients from $L$ such that $f(x_1,\dots,x_k)=0$. An infinite set $X$ of elements from $F$ is called  algebraically independent over $L$ if every finite subset of $X$ is algebraically independent over $L$. A maximal set of algebraically independent over $L$ elements from $F$ is called a transcendence basis of $F$ over $L$. The cardinality of a transcendence basis of $F$ over $L$ does not depend on this basis and is called the transcendence degree of $F$ over $L$. If $X$ is a transcendence basis of $F$ over $L$, then the subfield $L(X)$ of $F$ is isomorphic to the field of rational functions over the set of variables $X$ with coefficients from $L$. A field $F$ is said to be  algebraically closed if for every polynomial $f(t)$ of non-zero degree with coefficients from $F$ there exists an element $x$ in $F$ such that $f(x)=0$. The minimal algebraically closed field which contains $F$ is called the algebraic closure of $F$ and is denoted by $\overline{F}$. For every field there exists a unique (up to isomorphism) algebraic closure. If $F$ is an algebraically closed field with the prime subfield $L$ (i.~e. $L=F_p$ or $L=\mathbb{Q}$) and the trascendence basis $X$ of $F$ over $L$, then $F=\overline{L(X)}$, therefore every algebraically closed field is completely determined by its characteristic and the trascendence degree over the prime subfield. 

For $n\geq 2$ denote by ${\rm M}_n(F)$ the set of all $n\times n$ matrices over $F$. If $\varphi$ is an automorphism of $F$, then $\varphi$ induces a map $\varphi_n:{\rm M}_n(F)\to {\rm M}_n(F)$ which maps a matrix $A=(a_{i,j})\in {\rm M}_n(F)$ to $\varphi_n(A)=(\varphi(a_{i,j}))$. The map $\varphi_n$ induces an automorphism of classical linear groups (general linear group ${\rm GL}_n(F)$, special linear group ${\rm SL}_n(F)$, orthogonal group ${\rm O}_n(F)$, $\dots$) which we denote by the same symbol $\varphi_n$. Note that if an $n\times n$ matrix $A$ has rational coefficients, then $\varphi_n(A)=A$.

Let $F$ be an algebraically closed field of zero characteristic with countable transcendence degree over $\mathbb{Q}$. This field is isomorphic to $\overline{\mathbb{Q}(S)}$, where $S=\{s_1, s_2,\dots\}$ is a countable set of variables. In \cite[Theorem~6]{Nas4} it is proved that there exists an automorphism $\varphi$ of $F=\overline{\mathbb{Q}(S)}$ which induces an automorphism $\varphi_n$ of ${\rm GL}_n(F)$ with $R(\varphi_n)=1$. The following theorem gives some details about this automorphism $\varphi$ which follow from the proof of \cite[Theorem~6]{Nas4}.
\begin{theorem}\label{construction} Let $S$ be a countable set of variables, $F=\overline{\mathbb{Q}(S)}$, and $n\geq 1$ be a positive integer. Then there exists an automorphism $\varphi$ of $F$ (depending on $n$) such that for every matrix $A\in{\rm GL}_n(F)$ there exists an $n\times n$ matrix $X$ such that the entries of $X$ are algebraically independent over $\mathbb{Q}$ and $\varphi_n(X)=XA$.
\end{theorem}
Since entries of the matrix $X$ from Theorem~\ref{construction} are algebraically independent over $\mathbb{Q}$, and ${\rm det}(X)$ is the polynomial over the entries of $X$ with integer coefficients, we have ${\rm det}(X)\neq 0$ and $X\in {\rm GL}_n(F)$. So, from Theorem~\ref{construction} follows that for each matrix $A\in{\rm GL}_n(F)$ there exists a matrix $X\in{\rm GL}_n(F)$ such that $A=X^{-1}\varphi_n(X)$, i.~e. all matrices from  ${\rm GL}_n(F)$ are $\varphi_n$-conjugated to the identity matrix and $R(\varphi_n)=1$.
\section{Orthogonal groups over $\overline{\mathbb{Q}(s_1, s_2,\dots)}$}\label{sec3}
We use classical notation. Symbols $I_n$ and $O_{n\times m}$ denote the identity $n\times n$ matrix and the $n\times m$ matrix with zero entries, respectively. If  $A$ an $n\times n$ matrix and $B$ an $m\times m$ matrix over a field $F$, then the symbol  $A\oplus B$ denotes the direct sum of the matrices $A$ and $B$, i.~e. the block-diagonal $(m+n)\times (m+n)$ matrix
$$
\newcommand{\tempa}{\multicolumn{1}{c|}{A}}
\newcommand{\tempb}{\multicolumn{1}{|c}{B}}
\begin{pmatrix}
\tempa&O_{n\times m}\\\cline{1-2}
O_{m\times n}&\tempb
\end{pmatrix}
.
$$
It is clear that if $\varphi$ is an automorphism of $F$, then $\varphi_{n+m}(A\oplus B)=\varphi_n(A)\oplus\varphi_m(B)$. 
We denote by ${\rm GL}_n(F)$, ${\rm SL}_n(F)$, respectively, the group of $n\times n$ invertible matrices over $F$, and the group of $n\times n$ matrices over $F$ with determinant $1$.
The orthogonal group ${\rm O}_{n}(F)$ over  $F$ is the group
$${\rm O}_{n}(F)=\left\{A\in {\rm GL}_{n}(F) ~|~AA^T={\rm I}_n\right\},$$
where $^T$ denotes transpose. Denote by ${\rm SO}_n(F)={\rm O}_n(F)\cap {\rm SL}_n(F)$ the special orthogonal group, and by $\Omega_n(F)$ the derived subgroup of ${\rm O}_n(F)$. Symbols ${\rm PO}_n(F)$, ${\rm PSO}_n(F)$, ${\rm P}\Omega_n(F)$ denote, respectively, the quotients of the groups ${\rm O}_n(F)$, ${\rm SO}_n(F)$, $\Omega_n(F)$ by its centers. If $n\geq 2$, and $F$ is an algebraically closed field, then from \cite[Theorem 5.17]{Art} follows that $\Omega_n(F)={\rm SO}_n(F)$. The group of lower triangular matrices over $F$ is denoted by ${\rm T}_n(F)=\{(a_{i,j})\in {\rm GL}_{n}(F)~|~a_{i,i}\neq0,~\text{for}~i=1,\dots,n,~a_{i,j}=0~\text{for}~i<j\}$.

It is well known that every matrix $X$ from ${\rm GL}_n(\mathbb{R})$ can be presented as a product $X=LQ$ for some $L\in {\rm T}_n(\mathbb{R})$, $Q\in {\rm O}_n(\mathbb{R})$ (see, for example, \cite[Section~II, Lecture~7]{TreBau}). The following lemma says that some invertible matrices over an algebraically closed field of zero characteristic can be also presented in such form.
\begin{lemma}\label{gramm} Let $F$ be an algebraically closed field of zero characteristic, and $X$ be an $n\times n$ matrix such that entries of $X$ are algebraically independent over $\mathbb{Q}$. Then there exist matrices $L\in {\rm T}_n(F)$, $Q\in {\rm O}_n(F)$ such that $X=LQ$.
\end{lemma}
\noindent \textbf{Proof.} For $U=(u_1,\dots,u_n), V=(v_1,\dots,v_n)\in F^n$ denote by $\langle U,V\rangle$ the bilinear form $\langle U,V\rangle=u_1v_1+u_2v_2+\dots+u_nv_n$ on $F^n\times F^n$. Let $X_1,\dots,X_n$ be the rows of the matrix $X$. Denote by $Y_1,\dots,Y_n$ the following vectors from $F^n$.
\begin{align}
\notag Y_1&=X_1,\\
\notag Y_2&=X_2-\frac{\langle Y_1,X_2\rangle}{\langle Y_1, Y_1\rangle}Y_1\\
\label{gram} Y_3&=X_3-\frac{\langle Y_1,X_3\rangle}{\langle Y_1, Y_1\rangle}Y_1-\frac{\langle Y_2,X_3\rangle}{\langle Y_2, Y_2\rangle}Y_2\\
\notag &~~\vdots\\
\notag Y_{n}&=X_{n}-\frac{\langle Y_1,X_{n}\rangle}{\langle Y_1, Y_1\rangle}Y_1-\frac{\langle Y_2,X_{n}\rangle}{\langle Y_2, Y_2\rangle}Y_2-\dots-\frac{\langle Y_{n-1},X_{n}\rangle}{\langle Y_{n-1}, Y_{n-1}\rangle}Y_{n-1}.
\end{align}
The element $Y_1$ is uniquely defined in (\ref{gram}). In order the element $Y_2$ be well defined we need to check that  $\langle Y_1, Y_1\rangle\neq 0$. If $\langle Y_1, Y_1\rangle\neq 0$, then in order $Y_3$ be well defined we need to check that  $\langle Y_2, Y_2\rangle\neq 0$. So, in order formulas (\ref{gram}) be well defined, we need subsequently check that $\langle Y_k, Y_k\rangle\neq 0$ for all $k=1,\dots,n-1$. 

Since the entries of $X$ are algebraically independent over $\mathbb{Q}$, we can think about these entries as about variables, and we can think about $\langle Y_1, Y_1\rangle$ as about rational function from $\mathbb{Q}^{n\times n}$ to $\mathbb{Q}$ (elements in entries of $Y_1$ are rational functions in terms of entries of the matrix $X$, and $\langle Y_1, Y_1\rangle$ is a rational function over entries of $Y_1$). Denote this function by $f_1(X)=\langle Y_1, Y_1\rangle$. Since $f_1({\rm I}_n)=\langle(1,0,\dots,0), (1,0,\dots,0)\rangle=1$ , the function $f_1(X)$ is not constantly equal to zero, therefore $\langle Y_1, Y_1\rangle\neq 0$ as an element from $F$, and $Y_2$ is well defined. In a similar way, we can think about $\langle Y_2, Y_2\rangle$ as about rational function from $\mathbb{Q}^{n\times n}$ to $\mathbb{Q}$. Denote this function by $f_2(X)=\langle Y_2, Y_2\rangle$. Since $f_2({\rm I}_n)=\langle(0,1,0,\dots,0), (0,1,0,\dots,0)\rangle=1$, the function $f_2(X)$ is not constantly equal to zero, therefore $\langle Y_2, Y_2\rangle\neq 0$ as an element from $F$, and $Y_3$ is well defined. Using the same considerations, we prove that $\langle Y_k, Y_k\rangle\neq 0$ for all $k=1,\dots,n-1$ and therefore formulas (\ref{gram}) uniquely define $Y_1,\dots,Y_n$ from $X_1,\dots,X_n$.

For $1\leq i\neq j\leq n$ it is easy to check that $\langle Y_i, Y_j\rangle=0$ (formulas (\ref{gram}) are the same as in the classical Gramm-Schmidt orthogonalization process. See, for example, \cite[Section~II, Lecture~8]{TreBau}). So, denoting by $Q_k=Y_k/\sqrt{\langle Y_k,Y_k\rangle}$ (we can take square roots since $F$ is algebraically closed) we see that $\langle Q_i, Q_j\rangle=0$ if $i\neq j$ and $\langle Q_i, Q_i\rangle=1$ for $1\leq i, j\leq n$. It means that the matrix $Q$ with the first line $Q_1$, the second line $Q_2$, $\dots$, the $n$-th line $Q_n$, is orthogonal. From formulas (\ref{gram}) follows that $Q=TX$ for $T\in {\rm T}_n(F)$. Denoting by $L=T^{-1}$, we conclude the proof of the lemma.\hfill$\square$
\begin{lemma}\label{diagg} Let $S$ be a countable set of variables, $F=\overline{\mathbb{Q}(S)}$, $n\geq 1$ be an integer, and  $\varphi$ be an automorphism of $F$ introduced in Theorem~\ref{construction}. Then every diagonal matrix $D\in {\rm O}_n(F)$ is $\varphi_n$-conjugated in ${\rm O}_n(F)$ either to ${\rm I}_n$, or to ${\rm diag}(-1,\underbrace{1,\dots,1}_{n-1})$.
\end{lemma}
\noindent \textbf{Proof.} If $n=1$, then the statement is obvious. So, we can assume that $n\geq 2$. Since $D$ is orthogonal diagonal matrix, $D={\rm diag}(\varepsilon_1,\dots,\varepsilon_n)$ for $\varepsilon_1,\dots,\varepsilon_n\in\{\pm1\}$. Denote by $m$ the number of $-1$'s on the diagonal of $D$. Then there exists a permutation matrix $X$ such that 
$$XDX^{-1}={\rm diag}(\underbrace{-1,\dots,-1}_m,1,\dots,1)$$
By Theorem~\ref{construction} there exists a matrix $Y=(y_{i,j})$ such that the entries of $Y$ are algebraically independent over $\mathbb{Q}$ and $\varphi_n(Y)=-Y$. Since $n\geq 2$, there exist two algebraically independent over $\mathbb{Q}$ elements $\alpha=y_{1,1}$, $\beta=y_{1,2}$ such that $\varphi(\alpha)=-\alpha$, $\varphi(\beta)=-\beta$. Since $\alpha,\beta$ are algebraically independent over $\mathbb{Q}$, we have $\alpha^2+\beta^2\neq 0$. Denote by $Z$ the following $2\times 2$ matrix.
$$Z=\frac{1}{\sqrt{\alpha^2+\beta^2}}\begin{pmatrix}\alpha&\beta\\
-\beta&\alpha
\end{pmatrix}$$
If $m$ is even, then denote by 
$$T=\underbrace{Z\oplus\dots\oplus Z}_{m/2}\oplus{\rm I}_{n-m}.$$
Using direct calculations, it is clear that $T\in {\rm O}_n(F)$. From the formulas for matrix $Z$ follows that 
$$T^{-1}\varphi_n(T)={\rm diag}(\underbrace{-1,\dots,-1}_m,1,\dots,1)=XDX^{-1}.$$
Since $X$ is permutation matrix, $X\in {\rm O}_n(F)$ and $\varphi_n(X)=X$ (since entries of $X$ are integers). Therefore $D=(TX)^{-1}\varphi_n(TX)$, where $TX\in {\rm O}_n(F)$, i.~e. $D$ is $\varphi_n$-conjugated in ${\rm O}_n(F)$ to ${\rm I}_n$.

If $m$ is odd, then denote by 
$$T={\rm I}_1\oplus\underbrace{Z\oplus\dots\oplus Z}_{(m-1)/2}\oplus{\rm I}_{n-m}.$$
Using direct calculations, it is clear that $T\in {\rm O}_n(F)$. From the formulas for matrix $Z$ follows that 
$$T^{-1}{\rm diag}(-1,1,\dots,1)\varphi_n(T)={\rm diag}(\underbrace{-1,\dots,-1}_m,1,\dots,1)=XDX^{-1}.$$
Similar to the case when $m$ is even we see that $D=(TX)^{-1}{\rm diag}(-1,1,\dots,1)\varphi_n(TX)$, where $TX\in {\rm O}_n(F)$, i.~e. $D$ is $\varphi_n$-conjugated in ${\rm O}_n(F)$ to ${\rm diag}(-1,1,\dots,1)$.\hfill$\square$

\begin{theorem}\label{ort} Let $S$ be a countable set of variables, $F=\overline{\mathbb{Q}(S)}$, $n\geq 1$ be an integer,  and $\varphi$ be an automorphism of $F$ introduced in Theorem~\ref{construction}. Then $\varphi$ induces an automorphism $\varphi_n$ of ${\rm O}_n(F)$ with $R(\varphi_n)=2$.
\end{theorem}
\noindent \textbf{Proof.} At first, note that ${\rm I}_n$ and ${\rm diag}(-1,1,\dots,1)$ cannot be $\varphi_n$-conjugated. Indeed,  if ${\rm I}_n$ and ${\rm diag}(-1,1,\dots,1)$ are $\varphi_n$-conjugated, then there exists a matrix $X\in {\rm O}_n(F)$ such that  ${\rm diag}(-1,1,\dots,1)=X^{-1}\varphi_n(X)$. Therefore 
\begin{equation}\label{contr}-1={\rm det}({\rm diag}(-1,1,\dots,1))={\rm det}(X^{-1}){\rm det}(\varphi_n(X)).
\end{equation}
Since $X\in{\rm O}_n(F)$, ${\rm det}(X)=\varepsilon$, where $\varepsilon\in \{\pm 1\}$. Since $\varphi_n$ acts as $\varphi$ on entries of $X$, and ${\rm det}(X)$ is a polynomial with integer coefficients over entries of $X$, we have ${\rm det}(\varphi_n(X))=\varphi({\rm det}(X))=\varphi(\varepsilon)=\varepsilon$. From equality (\ref{contr}) follows that $-1=1$ and we have contradiction, i.~e. ${\rm I}_n$ and ${\rm diag}(-1,1,\dots,1)$ are not $\varphi_n$-conjugated and $R(\varphi_n)\geq 2$. Let us prove that $R(\varphi_n)=2$, i.~e. that every matrix $A\in{\rm O}_n(F)$ is $\varphi_n$-conjugated in ${\rm O}_n(F)$ either to ${\rm I}_n$ or to ${\rm diag}(-1,1,\dots,1)$. 

Let $A$ be a matrix from ${\rm O}_n(F)$. By Theorem~\ref{construction} there exists a matrix $X\in{\rm GL}_n(F)$ such that the entries of $X$ are algebraically independent over $\mathbb{Q}$ and $A=X^{-1}\varphi_n(X)$. From Lemma~\ref{gramm} follows that $X$ can be presented as a product $X=LQ$, where $L$ is a lower triangular matrix over $F$, and $Q$ is an orthogonal matrix over $F$. Since $A$ is orthogonal matrix, we have
\begin{align}
\notag {\rm I}_n&=AA^T\\
\notag &=X^{-1}\varphi_n(X)(X^{-1}\varphi_n(X))^T\\
\notag&=X^{-1}\varphi_n(X)\varphi_n(X)^T (X^{-1})^T\\
\notag&=(LQ)^{-1}\varphi_n(LQ)\varphi_n(LQ)^T ((LQ)^{-1})^T\\
\label{restr}&=Q^{-1}L^{-1}\varphi_n(L)\varphi_n(Q)\varphi_n(Q)^T\varphi_n(L)^T (L^{-1})^T(Q^{-1})^T
\end{align}
Since $Q$ is orthogonal, from equality (\ref{restr}) follows that ${\rm I}_n=L^{-1}\varphi_n(L)\varphi_n(L)^T (L^{-1})^T$ or
\begin{equation}\label{lrtr}
L^{-1}\varphi_n(L)=L^T\varphi_n(L^T)^{-1}.
\end{equation}
Since $L$ is a lower triangular matrix, the matrix on the left in equality (\ref{lrtr}) is lower triangular, and the matrix on the right in equality (\ref{lrtr}) is upper triangular. It is possible only in the case when $L^{-1}\varphi_n(L)=D$ is diagonal. Therefore we have the following equality.
\begin{equation}\label{diag}A=X^{-1}\varphi_n(X)=(LQ)^{-1}\varphi_n(L)\varphi(Q)=Q^{-1}D\varphi_n(Q).
\end{equation}
Since $A,Q$ are orthogonal, from equality (\ref{diag}) follows that $D$ is diagonal orthogonal matrix. From Lemma~\ref{diagg} follows that $D$ is $\varphi_n$-conjugated in ${\rm O}_n(F)$ either to ${\rm I}_n$ or to ${\rm diag}(-1,1,\dots,1)$. Therefore $A$ is $\varphi_n$-conjugated in ${\rm O}_n(F)$ either to ${\rm I}_n$ or to ${\rm diag}(-1,1,\dots,1)$.\hfill$\square$

\section{Symplectic groups over $\overline{\mathbb{Q}(s_1, s_2,\dots)}$}\label{sec4}
Denote by $J$, $\Omega_n$ the following $2\times 2$ and $2n\times 2n$ matrices 
\begin{align}\notag J=\begin{pmatrix}0&1\\-1&0\end{pmatrix},&&\Omega_n=\underbrace{J\oplus\dots\oplus J}_n.\end{align}
The symplectic group ${\rm Sp}_{2n}(F)$ over a field $F$ is the group
$${\rm Sp}_{2n}(F)=\left\{A\in {\rm GL}_{2n}(F) ~|~A\Omega_n A^T=\Omega_n\right\}.$$
It is easy to see that ${\rm Sp}_2(F)={\rm SL}_2(F)$. Denote by ${\rm PSp}_{2n}(F)$ the quotient of ${\rm Sp}_{2n}(F)$ by its center.

The set ${\rm M}_2(F)$ of all $2\times 2$ matrices over $F$ forms a (non commutative) ring and we can think about a matrix $A\in {\rm M}_{2n}(F)$ as about a matrix from ${\rm M}_n({\rm M}_2(F))$ over ${\rm M}_2(F)$ dividing $A$ into $2\times 2$ blocks. We say that a matrix $A\in {\rm M}_{2n}(F)={\rm M}_n({\rm M}_2(F))$ is a lower block triangular with $2\times 2$ blocks matrix if  it belongs to the set $\{(a_{i,j})\in {\rm M}_{n}({\rm M}_2(F)) ~|~a_{i,j}=O_{2\times2}~\text{for}~i<j\}$.
\begin{lemma}\label{gramm2} Let $F$ be an algebraically closed field of zero characteristic, and $X$ be a $2n\times 2n$ matrix such that the entries of $X$ are algebraically independent over $\mathbb{Q}$. Then there exist a matrix $Q\in {\rm Sp}_{2n}(F)$ and a lower block triangular with $2\times 2$ blocks matrix $L$ such that $X=LQ$.
\end{lemma}
\noindent \textbf{Proof.} For $U=(u_1,\dots,u_n), V=(v_1,\dots,v_n)\in {\rm M}_2(F)^n$ denote by $\langle U,V\rangle$ the following $2\times 2$ matrix with entries from $F$ 
$$\langle U,V\rangle=U\Omega_n V^T=u_1Jv_1^T+\dots+a_nJv_n^T.$$
We can think about $\langle U,V\rangle$ as about function from ${\rm M}_2(F)^n\times {\rm M}_2(F)^n\to {\rm M}_2(F)$. Using direct calculations it is easy to check that the equalities 
\begin{align}
\notag \langle U+V,W\rangle&=\langle U,W\rangle + \langle V,W
\rangle\\
\label{lin1} \langle xU,V\rangle&=x\langle U,V\rangle
\end{align}
hold for arbitrary $U, V, W\in {\rm M}_2(F)^n$, $x\in {\rm M}_2(F)$. For $U=(u_1,\dots,u_n)\in {\rm M}_2(F)^n$ denote by ${\rm d}(U)={\rm det}(u_1)+\dots+{\rm det}(u_n)$. From direct calculations follows that 
\begin{equation}\label{inv1}\langle U, U\rangle={\rm d}(U)J
\end{equation}

Let $X\in {\rm M}_{2n}(F)$ be a matrix from the formulation of the theorem. Express $X$ as a matrix over ${\rm M}_2(F)$ ($X\in {\rm M}_{2n}(F)={\rm M}_n({\rm M}_2(F))$) and let $X_1,\dots,X_n$ be the rows of the matrix $X$ (the entries of $X_1,\dots,X_n$ are $2\times 2$ matrices). Denote by $Y_1,\dots,Y_n$ the following elements from ${\rm M}_2(F)^n$.
\begin{align}
\notag Y_1&=X_1,\\
\notag Y_2&=X_2+{\rm d}(Y_1)^{-1}\langle X_2, Y_1\rangle J Y_1\\
\label{grams} Y_3&=X_3+{\rm d}(Y_1)^{-1}\langle X_3, Y_1\rangle J Y_1+{\rm d}(Y_2)^{-1}\langle X_3, Y_2\rangle J Y_2\\
\notag &~~\vdots\\
\notag Y_{n}&=X_{n}+{\rm d}(Y_1)^{-1}\langle X_n, Y_1\rangle J Y_1+\dots+{\rm d}(Y_{n-1})^{-1}\langle X_n, Y_{n-1}\rangle J Y_{n-1}.
\end{align}
The element $Y_1$ is uniquely defined in (\ref{grams}). In order the element $Y_2$ be well defined we need to check that  ${\rm d}(Y_1)\neq 0$. If ${\rm d}(Y_1)\neq 0$, then in order $Y_3$ be well defined we need to check that  ${\rm d}(Y_2)\neq 0$. So, in order formulas (\ref{grams}) be well defined, we need subsequently check that ${\rm d}(Y_k)\neq 0$ in $F$ for all $k=1,\dots,n-1$. 

Since the entries of $X$ are algebraically independent over $\mathbb{Q}$, we can think about these entries as about variables, and we can think about ${\rm d}(Y_1)$ as about rational function from $\mathbb{Q}^{2n\times 2n}$ to $\mathbb{Q}$. Denote this function by $f_1(X)={\rm d}(Y_1)$. Since 
$$f_1({\rm I}_{2n})={\rm d}\left(({\rm I}_2,O_{2\times 2},\dots,O_{2\times 2})\right)={\rm det}({\rm I}_2)=1,$$
 the function $f_1(X)$ is not constantly equal to zero, therefore ${\rm d}(Y_1)\neq 0$ as an element from $F$, and $Y_2$ is well defined. In a similar way, we can think about ${\rm d}(Y_2)$ as about rational function from $\mathbb{Q}^{2n\times 2n}$ to $\mathbb{Q}$. Denote this function by $f_2(X)={\rm d}(Y_2)$. Since $$f_2({\rm I}_{2n})=\left((O_{2\times 2}, {\rm I}_2, O_{2\times 2},\dots,O_{2\times 2})\right)={\rm det}({\rm I}_2)=1,$$ 
 the function $f_2(X)$ is not constantly equal to zero, therefore ${\rm d}(Y_2)\neq 0$ as an element from $F$, and $Y_3$ is well defined. In the same way we prove that ${\rm d}(Y_k)\neq 0$ for all $k=1,\dots,n-1$ and therefore formulas (\ref{gram}) uniquely define $Y_1,\dots,Y_n$ from $X_1,\dots,X_n$.

From equalities (\ref{lin1}), (\ref{inv1}) follows that 
\begin{align}
\notag \langle Y_2, Y_1\rangle&=\langle X_2+{\rm d}(Y_1)^{-1}\langle X_2, Y_1\rangle J Y_1, Y_1\rangle\\
\notag&\overset{(\ref{lin1})}{=}\langle X_2, Y_1\rangle+{\rm d}(Y_1)^{-1}\langle X_2, Y_1\rangle J\langle Y_1, Y_1\rangle\\
\notag&\overset{(\ref{inv1})}{=}\langle X_2, Y_1\rangle+{\rm d}(Y_1)^{-1}\langle X_2, Y_1\rangle J{\rm d}(Y_1)J\\
\label{ortogs}&=\langle X_2, Y_1\rangle-\langle X_2, Y_1\rangle=0.
\end{align}
Assuming that $\langle Y_i, Y_j\rangle=0$ for all $1\leq j<i$, using direct calculations similar to (\ref{ortogs}) we can show that $\langle Y_i, Y_j\rangle=0$ for all $1\leq j<i+1$. It means, that for all $1\leq j\neq i\leq n$ we have $\langle Y_i, Y_j\rangle=0$. 

Denote by $Y$ the matrix from ${\rm M}_n({\rm M}_2(F))$ with the first line $Y_1$, the second line $Y_2$, $\dots$, the $n$-th line $Y_n$. From equalities (\ref{grams}) follows that $Y=TX$, where $T$ is a lower block triangular with $2\times 2$ blocks matrix. For $k=1,\dots,n$ denote by  $Q_k=Y_k\sqrt{{\rm d}(Y_k)^{-1}}$ (we can take square roots since $F$ is algebraically closed). Since  $\langle Y_i, Y_j\rangle=0$ for all $1\leq j<i\leq n$, we have $\langle Q_i, Q_j\rangle=0$ for all $1\leq j<i\leq n$. Moreover, from equalities (\ref{lin1}), (\ref{inv1}) follows that $\langle Q_i, Q_i\rangle=J$ for all $i=1,\dots,n$. Therefore the matrix  $Q$ with the first line $Q_1$, the second line $Q_2$, $\dots$, the $n$-th line $Q_n$, belongs to ${\rm Sp}_{2n}(F)$. Since $Q_k=Y_k\sqrt{{\rm d}(Y_k)^{-1}}$ for $k=1,\dots,n$, we have $Q=DY=DTX$ for a diagonal $2n\times 2n$ matrix $D$. Denoting by $L=(DT)^{-1}$, we conclude the proof of the lemma.\hfill$\square$
\begin{lemma}\label{bdiag} Let $S$ be a countable set of variables, $F=\overline{\mathbb{Q}(S)}$, $n\geq 2$ be an even integer, and  $\varphi$ be an automorphism of $F$ introduced in Theorem~\ref{construction}. Let $D_1,\dots,D_{n/2}$ be matrices from ${\rm SL}_2(F)$ and $D=D_1\oplus\dots\oplus D_{n/2}$. Then $D$ is $\varphi_n$-conjugated in ${\rm Sp}_n(F)$ to ${\rm I}_n$.
\end{lemma}
\noindent \textbf{Proof.} For $k=1,\dots,n/2$ denote by $A_k=D_k\oplus {\rm I}_{n-2}$. By Theorem~\ref{construction} there exists a matrix $X_k\in{\rm GL}_n(F)$ such that the entries of $X_k$ are algebraically independent over $\mathbb{Q}$ and $\varphi_n(X_k)=X_kA_k$. Let
$$
X_k=\newcommand{\tempa}{\multicolumn{1}{c|}{P_k}}
\newcommand{\tempb}{\multicolumn{1}{|c}{S_k}}
\begin{pmatrix}
\tempa&Q_k\\\cline{1-2}
R_k&\tempb
\end{pmatrix}
,
$$
where $P_k$ is a $2\times 2$ matrix, $Q_k$ is a $2\times (n-2)$ matrix, $R_k$ is a $(n-2)\times 2$ matrix, and $S_k$ is $(n-2)\times (n-2)$ matrix. From the equality $\varphi_n(X_k)=X_kA_k$ and the fact that $A_k=D_k\oplus {\rm I}_{n-2}$ follows, that 
$\varphi_2(P_k)=P_kD_k$. Since the entries of $X_k$ are algebraically independent over $\mathbb{Q}$, and ${\rm det}(P_k)$ is a polynomial over the entries of $X$ with integer coefficients, ${\rm det}(P_k)\neq 0$ and $P_k$ is invertible. Therefore we have the equality
$$D_k=P_k^{-1}\varphi_2(P_k),$$
for $k=1,\dots,n/2$. From this equality follows that 
$$1={\rm det}(D_k)={\rm det}(P_k)^{-1}\varphi({\rm det}(P_k)).$$ Denote by $Y_k=\sqrt{{\rm det}(P_k)^{-1}}P_k$ (we can take square roots since $F$ is algebraically closed). Then ${\rm det}(Y_k)=1$ and 
$$Y_k^{-1}\varphi_2(Y_k)=\sqrt{{\rm det}(P_k)\varphi({\rm det}(P_k))^{-1}}P_k^{-1}\varphi_2(P_k)=D_k$$
for $k=1,\dots,n/2$. Denoting by $Y=Y_1\oplus\dots\oplus Y_{n/2}$ we have $D=Y^{-1}\varphi_n(Y)$, what concludes the proof.\hfill$\square$
\begin{theorem}\label{syms} Let $S$ be a countable set of variables, $F=\overline{\mathbb{Q}(S)}$, $n\geq 2$ be an even integer,  and $\varphi$ be an automorphism of $F$ introduced in Theorem~\ref{construction}. Then $\varphi$ induces an automorphism $\varphi_n$ of ${\rm Sp}_{n}(F)$ with $R(\varphi_n)=1$.
\end{theorem}
\noindent \textbf{Proof.} Let $A$ be a matrix from ${\rm Sp}_n(F)$. By Theorem~\ref{construction} there exists a matrix $X\in{\rm GL}_n(F)$ such that the entries of $X$ are algebraically independent over $\mathbb{Q}$ and $A=X^{-1}\varphi_n(X)$. From Lemma~\ref{gramm2} follows that $X$ can be presented as $X=LQ$, where $L$ is a lower block triangular with $2\times 2$ blocks matrix, and $Q\in{\rm Sp}_n(F)$. Since $A\in {\rm Sp}_n(F)$, we have
\begin{align}
\notag \Omega_{n/2}&=A\Omega_{n/2} A^T\\
\notag &=X^{-1}\varphi_n(X)\Omega_{n/2}(X^{-1}\varphi_n(X))^T\\
\notag&=X^{-1}\varphi_n(X)\Omega_{n/2}\varphi_n(X)^T (X^{-1})^T\\
\notag&=(LQ)^{-1}\varphi_n(LQ)\Omega_{n/2}\varphi_n(LQ)^T ((LQ)^{-1})^T\\
\label{restr1}&=Q^{-1}L^{-1}\varphi_n(L)\varphi_n(Q)\Omega_{n/2}\varphi_n(Q)^T\varphi_n(L)^T (L^{-1})^T(Q^{-1})^T
\end{align}
Since $Q\in {\rm Sp}_n(F)$, from equality (\ref{restr1}) follows that $\Omega_{n/2}=L^{-1}\varphi_n(L)\Omega_{n/2}\varphi_n(L)^T (L^{-1})^T$ or
\begin{equation}\label{lrtr1}
L^{-1}\varphi_n(L)=\Omega_{n/2}L^T\varphi_n(L^T)^{-1}\Omega_{n/2}^{-1}.
\end{equation}
Since $L$ is a lower block triangular with $2\times 2$ blocks matrix, and $\Omega_{n/2}$ is a block diagonal with $2\times 2$ blocks matrix, the matrix on the left in equality (\ref{lrtr1}) is a lower block triangular with $2\times 2$ blocks, and the matrix on the right in equality (\ref{lrtr1}) is an upper block triangular with $2\times 2$ blocks. It is possible only in the case when $L^{-1}\varphi_n(L)=D$ is block diagonal with $2\times 2$ blocks. Therefore we have the following equality.
\begin{equation}\label{diag1}A=X^{-1}\varphi_n(X)=(LQ)^{-1}\varphi_n(L)\varphi(Q)=Q^{-1}D\varphi_n(Q).
\end{equation}
Since $D$ is block diagonal with $2\times 2$ blocks, $D=D_1\oplus \dots\oplus D_{n/2}$  for $D_1,\dots,D_{n/2}\in {\rm M}_2(F)$. Since $A,Q\in{\rm Sp}_n(F)$, from equality (\ref{diag1}) follows that $D\in{\rm Sp}_n(F)$, therefore $D_k\in {\rm SL}_2(F)$ for all $k=1,\dots,n/2$. From Lemma~\ref{bdiag} follows that $D$ is $\varphi_n$-cojugated in ${\rm Sp}_{n}(F)$ to $I_n$, therefore $A$ is $\varphi_n$-conjugated in ${\rm Sp}_n(F)$  to ${\rm I}_n$.\hfill$\square$
\section{Orthogonal and symplectic groups in general case}\label{sec5}
Let us recall some necessary facts from model theory. The signature $\Sigma$ is a triple $(\mathcal{F},\mathcal{P},\rho)$, where $\mathcal{F}$ and $\mathcal{P}$ are disjoint sets not containing basic logical symbols, called, respectively, the set of function symbols and the set of predicate symbols, and $\rho:\mathcal{F}\cup \mathcal{P}\to \{0\}\cup \mathbb{N}$ is a function, called  arity, which assigns a non-negative integer  to every function or predicate symbol. A function symbol $f$ is called  $n$-ary if $\rho(f)=n$. A $0$-ary function symbol is called  a constant symbol.

The alphabet of a first-order logic of a signature $\Sigma=(\mathcal{F},\mathcal{P},\rho)$ consists of symbols of variables (usually $x$, $y$, $z$, $\dots$), logical operations (negation $\neg$, conjunction $\wedge$, disjunction $\vee$, implication $\to$), quantifiers (existential $\exists$, universal $\forall$), functional symbols from $\mathcal{F}$ (usually $f$, $g$, $\dots$ for symbols with positive arity, and $a$, $b$, $\dots$ for $0$-ary symbols), predicate symbols from $\mathcal{P}$ (usually $p$, $q$, $\dots$), parentheses, brackets and other punctuation symbols.

The set of terms of a first-order logic of a signature $\Sigma$ is inductively defined by the following rules: any variable is a term; for a functional symbol $f\in \mathcal{F}$ with $\rho(f)=n$ and terms $t_1$, $\dots$, $t_n$ the expression $f(t_1,\dots,t_n)$ is a term.

The set of formulas of a first-order logic of a signature $\Sigma$ is inductively defined by the following rules: for a predicate symbol $p\in \mathcal{P}$ with $\rho(p)=n$ and for terms $t_1$, $\dots$, $t_n$ the expression $p(t_1,\dots,t_n)$ is a formula; for terms $t_1$, $t_2$ the expression $t_1=t_2$ is a formula; if $x$ is a variable and $\varphi$, $\psi$ are formulas, then $\varphi\wedge\psi$, $\varphi\vee\psi$, $\varphi\to\psi$, $\neg \varphi$, $\exists x~\varphi$, $\forall x~\varphi$ are formulas. A set $\mathcal{A}$ of some formulas of a first-order logic of a signature $\Sigma$ is called  a theory of a signature $\Sigma$.

Let $\Sigma=(\mathcal{F}, \mathcal{P}, \rho)$ be a signature, $M$ be a non-empty set and $\sigma$ be a function which maps every functional symbol $f$ from $\mathcal{F}$ with $\rho(f)=n$ to $n$-ary function $\sigma(f):M^n\to M$, and maps every predicate symbol $p$ from $\mathcal{P}$ to $n$-ary relation $\sigma(p)\subseteq M^n$. Denote by $\mathcal{M}$ the pair $(M,\sigma)$. Let $s$ be a function which maps any variable to some element from $M$.  The interpretation $[[t]]_s$ of a term $t$ in $M$ with respect to $s$ is inductively defined by the following rules: $[[x]]_s=s(x)$ if $x$ is a variable and $[[f(t_1,\dots,t_n)]]_s=\sigma(f)([[t_1]]_s,\dots,[[t_n]]_s)$ for a functional symbol $f\in \mathcal{F}$ with $\rho(f)=n$ and terms $t_1$, $\dots$, $t_n$.  The truth of a formula $\varphi$ in $\mathcal{M}$ with respect to $s$ (we write $\mathcal{M}\models_s\varphi$ if $\varphi$ is true in $\mathcal{M}$ with respect to $s$) is inductively defined by the following rules:
 \begin{itemize}
 \item $\mathcal{M}\models_s p(t_1,\dots,t_n)$ if an only if $\left([[t_1]]_s,\dots,[[t_n]]_s\right)\in\sigma(p)$,
  \item $\mathcal{M}\models_s t_1=t_2$ if and only if $[[t_1]]_s=[[t_2]]_s$,
   \item $\mathcal{M}\models_s\varphi\wedge\psi$ if and only if $\mathcal{M}\models_s\varphi$ and $\mathcal{M}\models_s\psi$,
   \item $\mathcal{M}\models_s\varphi\vee\psi$ if and only if $\mathcal{M}\models_s\varphi$ or $\mathcal{M}\models_s\psi$,
    \item $\mathcal{M}\models_s\varphi\to\psi$ if and only if $\mathcal{M}\models_s\varphi$ implies $\mathcal{M}\models_s\psi$,
     \item $\mathcal{M}\models_s\neg \varphi$ if and only if $\mathcal{M}\models_s\varphi$ is not true,
      \item $\mathcal{M}\models_s\exists x~\varphi$ if and only if $\mathcal{M}\models_{s^{\prime}}\varphi$ for some $s^{\prime}$ with $s^{\prime}(y)=s(y)$ for all $y\neq x$,
       \item $\mathcal{M}\models_s\forall x~\varphi$ if and only if $\mathcal{M}\models_{s^{\prime}}\varphi$ for all $s^{\prime}$ with $s^{\prime}(y)=s(y)$ for all $y\neq x$.
 \end{itemize}
We say that $\mathcal{M}=(M,\sigma)$ is a model for a theory $\mathcal{A}$ if $\mathcal{M}\models_s\varphi$ for all formulas $\varphi\in \mathcal{A}$ and all functions $s$. A theory can have no models, one model or several models (even an infinite number). 
\begin{theorem}[L\"{o}wenheim-Skolem Theorem]  If a countable theory $\mathcal{A}$ of a signature $\Sigma$ has an infinite model, then for every infinite cardinal number $\kappa$ it has a model $\mathcal{M}=(M,\sigma)$ with $|M|=\kappa$.
\end{theorem}
\begin{theorem}\label{mthO}Let $F$ be an algebraically closed field of zero characteristic with infinite transcendence degree over $\mathbb{Q}$. Then there exists an automorphism $\varphi$ of $F$ which induces an automorphism $\varphi_n$ of ${\rm O}_n(F)$ with $R(\varphi_n)=2$. 
\end{theorem}
\noindent\textbf{Proof.} If the transcendence degree of $F$ over $\mathbb{Q}$ is countable, then $F=\overline{\mathbb{Q}(S)}$, where $S$ is a countable set of variables, and the result follows from Theorem~\ref{ort}. So let the transcendence degree of $F$ over $\mathbb{Q}$ be uncountable.

Let $\Sigma=(\mathcal{F},\mathcal{P},\rho)$ be a signature with $\mathcal{\mathcal{F}}=\{+,\cdot,-,^{-1},f,0,1\}$, $\mathcal{P}=\varnothing$ and $\rho(+)=\rho(\cdot)=2$, $\rho(-)=\rho(^{-1})=\rho(f)=1$, $\rho(0)=\rho(1)=0$. For terms $t_1$, $t_2$ of the signature $\Sigma$ we will write $t_1+t_2$, $t_1\cdot t_2$, $-t_1$, $t_1^{-1}$ instead of $+(t_1,t_2)$, $\cdot(t_1, t_2)$, $-(t_1)$, $^{-1}(t_1)$, respectively, and we will write $t_1\neq t_2$ instead of $\neg(t_1=t_2)$. Let $\mathcal{A}$ be the theory consisting of the following formulas.
\begin{enumerate}
\item $\forall x~\forall y~\forall z~[x+y=y+x]\wedge[x+(y+z)=(x+y)+z]$,
\item $\forall x~\forall y~\forall z~[x\cdot y=y\cdot x]\wedge[x\cdot(y\cdot z)=(x\cdot y)\cdot z]$,
\item $\forall x~\forall y~\forall z~x\cdot(y+z)=x\cdot y+x\cdot z$,
\item $\forall x~[x+0=x]\wedge[x\cdot1=x]$,
\item $\forall x~[x+(-x)=0]\wedge[(x\neq0)\to (x\cdot x^{-1}=1)]$,
\item $1\neq0$, $1+1\neq0$, $1+1+1\neq0$, $\dots$
\item $\forall y_1~\forall y_0~[(y_1\neq0)\to(\exists x~y_1\cdot x+y_0=0)]$,\\
$\forall y_2~\forall y_1~\forall y_0~[(y_2\neq0)\to(\exists x~y_2\cdot x\cdot x+y_1\cdot x+y_0=0)]$,\\
 $\dots$
\item $\forall x~\forall y~\exists z~[(f(x)=f(y))\to(x=y)]\wedge[f(z)=x]$,
\item $\forall x~\forall~y~[f(x+y)=f(x)+f(y)]\wedge[f(x\cdot y)=f(x)\cdot f(y)]$,
\item $\forall x_{1,1}~\forall x_{1,2}\dots\forall x_{1,n}~\forall x_{2,1}\dots \forall x_{2,n}\dots\forall x_{n,n}\\
\exists y_{1,1}~\exists y_{1,2}\dots\exists y_{1,n}~\exists y_{2,1}\dots \exists y_{2,n}\dots\exists y_{n,n}$\\
$[XX^T=I_{n}]\to[(YY^T=I_n)\wedge((f(Y)=YX)\vee(((-1)\oplus {\rm I_{n-1}}))f(Y)=YX))]$.\\
where $X=(x_{i,j})$, $Y=(y_{i,j})$ and $f(Y)=(f(y_{i,j}))$.
\end{enumerate}
Note that in formula 10 the expression $XX^T={\rm I}_n$ (multiplication of matrices) meas the conjunction of $n^{2}$ formulas (in terms of entries of matrices), and each of these $n^2$ formulas can be uniquely written in terms of $\cdot$, $+$ and $-$ (we write $XX^T={\rm I}_n$ for the simplicity of denotation). Similarly, formulas $YY^T=I_n$, $f(Y)=YX$, $((-1)\oplus {\rm I}_{n-1})f(Y)=YX$ can be written as conjunctions of $n^2$ formulas. So, what is written in 10 is a formula.

If there exists a model $\mathcal{M}=(M,\sigma)$ for $\mathcal{A}$, then formulas 1 and 2 say that addition and multiplication are commutative and associative,  formula 3 describes distributive law between addition and multiplication, formula 4 says that $0$ and $1$ are neutral elements with respect to addition and multiplication respectively, formula 5 says that $-x$ is an opposite to $x$ element and that $x^{-1}$ is an inverse to non-zero element $x$ element. So, all together formulas 1-5 say that $(M,\sigma(+),\sigma(\cdot))$ is a field with the zero element $\sigma(0)$ and the unit element $\sigma(1)$. The countable set of formulas 6 says that the field $M$ has zero characteristic. The countable set of formulas 7 says that $M$ is an algebraically closed field. Formula 8 says that $\sigma(f)$ is a bijection on $M$ and formula 9 says that this bijection respects addition and multiplication. So, formulas 8 and 9 together say that $\sigma(f)$ is an automorphism of the field $M$. Finally, formula 10 says that for every orthogonal matrix $X=(x_{i,j})$ with coefficients from $M$ there exists an orthogonal matrix $Y=(y_{i,j})$ with coefficients from $M$ such that either $(\sigma(f)(y_{i,j}))=YX$ or $((-1)\oplus {\rm I}_{n-1})(\sigma(f)(y_{i,j}))=YX$. In other words all together formulas 1-10 say that $M$ is an algebraically closed field of zero characteristic and $\sigma(f)$ is an automorphism of $M$ which induced an automorphism $\varphi$  of ${\rm O}_n(M)$ with $R(\varphi)=2$ (every orthogonal matrix is $\varphi$-conjugated either to ${\rm I}_n$ or to $(-1)\oplus {\rm I}_{n-1}={\rm diag}(-1,1,\dots,1)$).

The previous paragraph is written in assumption that the theory $\mathcal{A}$ has a model, but by Theorem \ref{ort} it has an infinite model. Since $\mathcal{A}$ is a countable theory, by L\"{o}wenheim-Skolem theorem it has a model $\mathcal{M}=(M,\sigma)$ with $|M|=\kappa=|F|$, where $F$ is a field from the formulation of the theorem. Since $|M|=\kappa$ is uncountable, the transcendence degree of $M$ over $\mathbb{Q}$ is infinite and is equal to $\kappa$. Since every algebraically closed field is completely determined by its characteristic and transcendence degree over prime subfield, we have $M\cong F$. Therefore there exists an automorphism $\varphi$ of the field $F$ which induces an automorphism of ${\rm O}_n(F)$ with Reidemeister number equals to $2$.\hfill$\square$
\begin{cor}\label{mthSO}Let $F$ be an algebraically closed field of zero characteristic with infinite transcendence degree over $\mathbb{Q}$. Then there exists an automorphism $\varphi$ of $F$ which induces an automorphism $\varphi_n$ of ${\rm SO}_n(F)$ with $R(\varphi_n)=1$. 
\end{cor}
\noindent\textbf{Proof.} Let $\varphi$ be an automorphism from Theorem~\ref{mthO}, and $A$ be an arbitrary matrix from ${\rm SO}_n(F)$. From the proof of Theorem~\ref{mthSO} follows that $A$ is $\varphi_n$ conjugated in ${\rm O}_n(F)$ either to ${\rm I}_n$ or to ${\rm diag}(-1,1,\dots,1)$. Since ${\rm det}(A)=1$, the matrix $A$ cannot be $\varphi_n$-conjugated in ${\rm O}_n(F)$ to ${\rm diag}(-1,1,\dots,1)$, therefore $A$ is $\varphi_n$-conjugated in ${\rm O}_n(F)$ to ${\rm I}_n$, i.~e. there exists a matrix $X\in {\rm O}_n(F)$ such that $A=X^{-1}\varphi_n(X)$. Denote by $Y=(({\rm det}(X))\oplus {\rm I}_{n-1})X$. Then $Y\in{\rm SO}_n(F)$ and $Y^{-1}\varphi_n(Y)=X^{-1}\varphi_n(X)=A$, i.~e. an arbitrary matrix $A\in {\rm SO}_n(F)$ is $\varphi_n$-conjugated in ${\rm SO}_n(F)$ to ${\rm I}_{n}$.\hfill$\square$

Slightly modifying the set of formulas in the proof of Theorem~\ref{mthO} we get the following result.
\begin{theorem}\label{mthS}Let $F$ be an algebraically closed field of zero characteristic with infinite transcendence degree over $\mathbb{Q}$, and $n$ be an even integer. Then there exists an automorphism $\varphi$ of $F$ which induces an automorphism $\varphi_n$ of ${\rm Sp}_n(F)$ with $R(\varphi_n)=1$. 
\end{theorem}
\noindent\textbf{Proof.} The proof is the same as the proof of Theorem~\ref{mthO} changing formula 10 from the proof of Theorem~\ref{mthO} to the formula
\begin{align}
\notag&\forall x_{1,1}~\forall x_{1,2}\dots\forall x_{1,n}~\forall x_{2,1}\dots \forall x_{2,n}\dots\forall x_{n,n},
\exists y_{1,1}~\exists y_{1,2}\dots\exists y_{1,n}~\exists y_{2,1}\dots \exists y_{2,n}\dots\exists y_{n,n}\\
\notag&[X\Omega_{n/2}X^T=\Omega_{n/2}]\to[(Y\Omega_{n/2}Y^T=\Omega_{n/2})\wedge(f(Y)=YX)],
\end{align}
where $X=(x_{i,j})$, $Y=(y_{i,j})$ and $f(Y)=(f(y_{i,j}))$ and using Theorem~\ref{syms} instead of Theorem~\ref{ort} to provide an infinite model which we need in order to use L\"{o}wenheim-Skolem theorem.\hfill$\square$
\section{Chevalley groups}\label{sec6}
The following statement is almost obvious, it can be found it  \cite[Lemmas 2.1]{MubSan2}.
 \begin{lemma}\label{npr3} Let $G$ be a group, $\varphi$ be an automorphism of $G$ and $N$ be a $\varphi$ admissible normal subgroup of $G$. Denote by $\overline{\varphi}$ an automorphism of $G/N$ induced by $\varphi$. Then $R(\overline{\varphi})\leq R(\varphi)$.
\end{lemma}
\begin{theorem}\label{main} Let $F$ be an algebraically closed field of zero characteristic which has infinite transcendence degree over $\mathbb{Q}$. Then Chevalley groups of types $B_n, C_n, D_n$ over $F$ do not possess the $R_{\infty}$-property.
\end{theorem}
\noindent\textbf{Proof.} We separately consider the root systems $B_n, D_n$ and the root systems $C_n$.

\textbf{Case 1:} The root system has the type $B_n$ or $D_n$. Let $G$ be a Chevalley group of type $B_n$  or $D_n$  over $F$, and let $Z(G)$ be the center of $G$ (it is known to be finite).   The quotient group $G/Z(G)$ is isomorphic to ${\rm P}\Omega_{k}(F)={\rm PSO}_k(F)$, where $k=2n+1$ if the root system has the type $B_n$, and $k=2n$ if the root system has the type $D_n$ (see  \cite[\S 11.3, \S12.1]{Car}). Denote by $\varphi$ the automorphism of $F$ from Corollary~\ref{mthSO} (changing $n$ by $k$ in  Corollary~\ref{mthSO}). Denote by $\psi$, $\overline{\psi}$, $\theta$, $\overline{\theta}$, respectively, the automorphisms of $G$, $G/Z(G)$, ${\rm SO}_{k}(F)$, ${\rm PSO}_{k}(F)$ induced by $\varphi$. Since $G/Z(G)={\rm PSO}_{k}(F)$, we have $R(\overline{\theta})=R(\overline{\psi})$.  From Corollary~\ref{mthSO} we have $R(\theta)=1$. From Lemma~\ref{npr3} follows that $R(\overline{\theta})=1$, therefore $R(\overline{\psi})=1$, i.~e. every element from $G/Z(G)$ is $\overline{\psi}$-conjugated to the unit element in $G/Z(G)$. Hence every element from $G$ is $\psi$-conjugated to some element from $Z(G)$. Since $Z(G)$ is finite, the number of $\psi$-conjugacy classes in $G$ is finite, i.~e. $G$ does not possess the $R_{\infty}$-property.

\textbf{Case 2:} The root system has the type $C_n$. Let $G$ be a Chevalley group of type $C_n$ over $F$, and let $Z(G)$ be the center of $G$ (it is known to be finite).   The quotient group $G/Z(G)$ is isomorphic to ${\rm PSp}_{2n}(F)$ (see  \cite[\S 11.3, \S12.1]{Car}). Denote by $\varphi$ the automorphism of $F$ from Theorem~\ref{mthS} (changing $n$ by $2n$ in  Theorem~\ref{mthS}). Denote by $\psi$, $\overline{\psi}$, $\theta$, $\overline{\theta}$, respectively, the automorphisms of $G$, $G/Z(G)$, ${\rm Sp}_{2n}(F)$, ${\rm PSp}_{2n}(F)$ induced by $\varphi$. Since $G/Z(G)={\rm PSp}_{2n}(F)$, we have $R(\overline{\theta})=R(\overline{\psi})$.  From Theorem~\ref{mthS} we have $R(\theta)=1$. From Lemma~\ref{npr3} follows that $R(\overline{\theta})=1$, therefore $R(\overline{\psi})=1$, i.~e. every element from $G/Z(G)$ is $\overline{\psi}$-conjugated to the unit element in $G/Z(G)$. Hence every element from $G$ is $\psi$-conjugated to some element from $Z(G)$. Since $Z(G)$ is finite, the number of $\psi$-conjugacy classes in $G$ is finite, i.~e. $G$ does not possess the $R_{\infty}$-property.\hfill$\square$
 
\begin{theorem}\label{main2} Let $F$ be an algebraically closed field of zero characteristic. Then Chevalley groups of types $A_n,B_n, C_n, D_n$ over $F$ possess the $R_{\infty}$-property if and only if the transcendence degree of $F$ over $\mathbb{Q}$ is finite.
\end{theorem}
\noindent\textbf{Proof.} The necessity condition follows from Theorem~\ref{main} and \cite[Theorem~7]{Nas4}. The sufficiency condition follows from \cite[Theorem~3.2]{FelNas}.\hfill$\square$

~\\

~\\

{\small
\begin{spacing}{0.5}

\end{spacing}

~\\
KU Leuven KULAK, Etienne Sabbelaan 53, 8500 Kortrijk, Belgium, timur.nasybullov@mail.ru}
\end{document}